\documentclass{article}

    \usepackage[preprint]{neurips_2022}

\usepackage[utf8]{inputenc} 
\usepackage[T1]{fontenc}    
\usepackage{hyperref}       
\usepackage{url}            
\usepackage{booktabs}       
\usepackage{amsfonts}       
\usepackage{nicefrac}       
\usepackage{microtype}      
\usepackage{xcolor}         

\usepackage{bm}
\usepackage{amsmath}
\usepackage{amsthm}
\usepackage{amssymb}
\usepackage{algorithm}
\usepackage{algorithmic}
\usepackage{subcaption}
\usepackage{mathtools}
\usepackage{booktabs}
\usepackage{bbm}
\usepackage{courier}
\usepackage{adjustbox}
\newcommand\numberthis{\addtocounter{equation}{1}\tag{\theequation}}
\newcommand\myeq{\stackrel{\mathclap{\normalfont\mbox{\tiny i.i.d.}}}{=}}

\theoremstyle{definition}

\newtheorem{theorem}{Theorem}[section]

\newtheorem{lemma}[theorem]{Lemma}

\newtheorem*{remark}{Remark}

\title{Optimizing Objective Functions from Trained ReLU Neural Networks via Sampling}

\author{%
  Georgia Perakis \\
   Operations Research Center \\
  Massachusetts Institute of Technology \\
  Cambridge, MA 02142\\
   \texttt{georgiap@mit.edu} \\
   \And
  Asterios Tsiourvas \\
   Operations Research Center \\
  Massachusetts Institute of Technology \\
  Cambridge, MA 02142\\
   \texttt{atsiour@mit.edu} \\
}

\begin{document}

\maketitle

\begin{abstract}This paper introduces scalable, sampling-based algorithms that optimize trained neural networks with ReLU activations. We first propose an iterative algorithm that takes advantage of the piecewise linear structure of ReLU neural networks and reduces the initial mixed-integer optimization problem (MIP) into multiple easy-to-solve linear optimization problems (LPs) through sampling. Subsequently, we extend this approach by searching around the neighborhood of the LP solution computed at each iteration. This scheme allows us to devise a second, enhanced algorithm that reduces the initial MIP problem into smaller, easier-to-solve MIPs. We analytically show the convergence of the methods and we provide a sample complexity guarantee. We also validate the performance of our algorithms by comparing them against state-of-the-art MIP-based methods. Finally, we show computationally how the sampling algorithms can be used effectively to warm-start MIP-based methods.
\end{abstract}

\section{Introduction} \label{sec:introduction}

\par In recent years, deep learning has proven to be an extremely powerful tool for solving a wide variety of problems in many areas ranging from image recognition, and natural language processing to robotics among many others (\cite{lecun2015deep, deng2014deep, pouyanfar2018survey}). Even though many applications have mostly focused on the predictive part, optimizing over trained neural networks tractably is a question that has recently emerged. In this case, the end goal is to solve an optimization problem that maximizes an objective function from a trained neural network. Deep reinforcement learning problems with a large high-dimensional action space are an important application of such an approach. There the cost-to-go function or the state transition functions can be learned through neural networks (\cite{arulkumaran2017deep}). Therefore, finding the optimal decision requires solving an optimization problem with a neural network on the objective. Overall, optimizing trained neural networks is a hard task due to their highly nonlinear and non-convex nature. 

In what follows, we present the underlying optimization problem. In particular, we represent a neural network as a function $f(\bm{x})$ with decision variables $\bm{x}$ belonging to a bounded polyhedron $\mathcal{P}$. This gives rise to the following optimization problem:
\begin{equation}
    \begin{array}{ll@{}ll}
    \max\limits_{\bm{x}} &  f(\bm{x}) \\
    \text{subject to} & \bm{x} \in \mathcal{P}.
    \end{array} \label{eq:init}
\end{equation}

In this work, we focus on neural networks with ReLU activation functions ($\sigma(x) = \max\{x,0\}$).  ReLU activations are among the most widely used activations in deep learning architectures (\cite{agarap2018deep}) and also, allow us to express problem \eqref{eq:init} as a MIP. Recent research (\cite{fischetti2018deep,anderson2020strong,tsay2021partition}) has focused on improving existing formulations and on implementing such MIP approaches for a wide variety of applications, such as verifying the robustness of an output (\cite{tjeng2017evaluating,rob1}), compressing neural networks (\cite{serra2020lossless}) and counting linear regions of a network (\cite{serra2018bounding}) among others. Nevertheless, one disadvantage of MIP-based approaches is that they do not necessarily scale easily as the size of the network increases. In what follows, we present the contributions of this paper.

\subsection{Contributions}

{\bf Sampling-based algorithm:} In Section~\ref{sec:sampling_method}, we propose a novel, sampling-based, iterative algorithm for optimizing objective functions that come from trained ReLU neural networks in a computationally tractable way. The algorithm takes advantage of the piecewise linear structure of the output of the network and reduces the initial MIP into multiple LPs through sampling. 

 {\bf An enhanced algorithm that uses fewer samples:} In Section~\ref{sec:local}, we extend the proposed algorithm by searching for improved solutions in the neighborhood of the LP solution computed at each iteration. This approach trades off the number of LPs, which is equal to the number of samples, the originally proposed algorithm needs to solve with fewer and smaller --- in terms of binary variables --- MIPs. This enhanced approach improves the performance of the original algorithm and requires fewer samples.

 {\bf Analytical guarantees:} In Section~\ref{sec:proovable_guarantees}, we prove analytically the algorithms' convergence to the optimal solution and establish a sample complexity guarantee. We also show that the enhanced sampling algorithm finds a local optimum from its first iteration.
 
 \textbf{Computational tractability}: In Section~\ref{sec:comparison}, we evaluate computationally the proposed algorithms relative to the improved MIP formulation, as described by \cite{fischetti2018deep}, the big-$M$ + cuts approach proposed by \cite{anderson2020strong} and the $N=\{2,4\}$ equal-size partition by \cite{tsay2021partition}. Sampling-based algorithms are competitive in terms of performance, with the enhanced algorithm consistently retrieving the highest objective value across all methods. Furthermore, in Section~\ref{sec:warm} we show that the sampling algorithms can be used effectively to warm-start the aforementioned MIP-based methods.

\subsection{Related Work}

Recent research has focused on improving existing MIP formulations for optimizing ReLU neural networks. \cite{fischetti2018deep} modeled neural networks with maximum-based nonlinear operators as a MIP and described a bound-tightening technique to ease the solution of the formulation. \cite{anderson2020strong} presented a general framework that provides a way to construct ideal or sharp formulations for optimizing linear functions over general polyhedra and applied it to derive strong MIP formulations for ReLU neural networks. A related dual algorithm was proposed by \cite{de2021scaling}. \cite{tsay2021partition} introduced partition-based MIP formulations for optimizing ReLU neural networks that balance the model size and the tightness of the MIP formulation using ideas from disjunctive programming. Other methods for obtaining strong relaxations of ReLU neural networks include LP (\cite{lp1,lp2,lp3}), semidefinite programming (\cite{sem1,sem2}), and Langrangian decomposition (\cite{la1}) among others. Our work departs from the aforementioned approaches since it does not focus on improving MIP formulations, but uses them to derive tractable algorithms that optimize ReLU neural networks by reducing the initial MIP either into many easy-to-solve LPs or a few small MIPs via sampling. Furthermore, our work also shows how these algorithms can be used to warm-start MIP-based methods.

\section{Mixed-Integer Optimization Formulation}\label{sec:MIP}
In this paper, for ease of exposition, we focus on proposing methods that optimize the feed-forward neural network architecture (\cite{svozil1997introduction}) with one-dimensional continuous output and ReLU activations. Nevertheless, the methods we propose can be incorporated into more complex architectures, as long as they utilize ReLU activations (or more generally maximum-based nonlinear activations). In what follows, we describe the MIP approach (\cite{fischetti2018deep,serra2018bounding}) for problem \eqref{eq:init}. 

\subsection{One Layer} \label{sec:one_layer}
Let $f(\bm{x})$ to be the output of a one-layer feed-forward ReLU neural network with $n$ neurons under the $d$-dimensional input vector $x \in \mathcal{P} \subset \mathbb{R}^d$, where $\mathcal{P}$ is a bounded polyhedron. We describe $f(\bm{x})$ as
\begin{align*}
    f(\bm{x}) &= \sum_{i=1}^n a_i ReLU(\bm{w}_i^T \bm{x} + b_i), \label{eq:1_layer_relu} \numberthis
\end{align*}
\noindent where $\bm{w}_i \in \mathbb{R}^d, b_i \in \mathbb{R}, a_i \in \mathbb{R},  \ \forall i \in [n]$, are the weights of the network.

\par To obtain a MIP formulation for problem \eqref{eq:init}, with $f$ being a one-layer ReLU network, we rewrite the ReLU activation as the following set of mixed-integer linear constraints:
\begin{align*}
    y_i &=ReLU(\bm{w}_i^T\bm{x}+b_i) \iff y_i \in \mathcal{C}(\bm{x},\bm{w}_i,b_i), \numberthis
\end{align*}
\noindent where 
\begin{align}
    \mathcal{C}(\bm{x},\bm{w}_i,b_i) = \{y| \ y\geq 0, \  y\geq \bm{w}_i^T\bm{x} +b_i, \ y \leq u_iz_i, \ y \leq \bm{w}_i^T\bm{x} +b_i - l_i(1-z_i), \ z_i \in \{0,1\} \}.
\end{align}
In the previous definition, $z_i$ is a binary variable that indicates whether the ReLU activation of neuron $i$ is activated. Furthermore, we have that $u_i = \sup\limits_{\bm{x} \in \mathcal{P}} \{ \bm{w}_i^T\bm{x} + b_i \}$ and $l_i = \inf\limits_{\bm{x} \in \mathcal{P}} \{ \bm{w}_i^T\bm{x} + b_i \} $  are the maximum and the minimum values that the output of neuron $i$ can take. Both $u_i$ and $l_i$ are bounded when $\mathcal{P}$ is bounded and they can be calculated efficiently using linear optimization. 

Using the previous results, problem \eqref{eq:init} can be reformulated as the following MIP:
\begin{align*}
    \max\limits_{\bm{x}\in\mathcal{P}} \sum_{i=1}^n a_i y_i, \ s.t.  \  y_i \in \mathcal{C}(\bm{x},\bm{w}_i, b_i), \ \forall i \in [n]. \numberthis \label{eq:1_layer_mip}
\end{align*}

\subsection{Multiple Layers}

Let $f(\bm{x})$ be a $K$-layer feed-forward ReLU neural network. We denote the $d$-dimensional input vector $\bm{x}$ as $\bm{x}_1 $ for ease of notation. The output of the network can be written as:
\begin{align*}
    f(\bm{x}_1) &= \sum_{i=1}^{n_K} a_{i} \bm{x}_{K,i}, \ \bm{x}_k = ReLU(\bm{W}_{k-1} \bm{x}_{k-1} + \bm{b}_{k-1}), \ k =2,\dots,K, \numberthis
\end{align*}
\noindent where $\bm{W}_{k} \in \mathbb{R}^{n_{k+1}\times n_{k}}, \bm{b}_{k} \in \mathbb{R}^{n_{k+1}}, \ \forall k=1,\dots,K-1,$ and $  \bm{a} \in \mathbb{R}^{n_K} $ are the weights of the network, with $n_k$ being the number of neurons at layer $k$ (we assume $n_1=d$). Therefore, problem \eqref{eq:init} becomes:
\begin{equation}
    \begin{array}{ll@{}ll}
    \max\limits_{\bm{x} } &  \sum_{i=1}^{n_K} a_{i} ReLU((\bm{W}_{K-1}\bm{x}_{K-1})_i + b_{K-1,i}) \\
    \text{subject to} & \bm{x}_k =  ReLU(\bm{W}_{k-1} \bm{x}_{k-1} + \bm{b}_{k-1}), \ k =2,\dots,K,\\ 
    & \bm{x}_1 \in \mathcal{P}.
    \end{array}
\end{equation}

By formulating ReLU activations as mixed-integer linear constraints, as described in Section~\ref{sec:one_layer}, problem \eqref{eq:init} can be reformulated as the following MIP:
\begin{equation}
    \begin{array}{ll@{}ll}
    \max\limits_{\bm{x}, \bm{z}} &  \sum_{i=1}^{n_K} a_{i} x_{K,i} \\
    \text{subject to} & \bm{x}_{k} \leq \bm{W}_{k-1} \bm{x}_{k-1} + \bm{b}_{k-1} -\bm{l}_{k}\odot(1-\bm{z}_{k}), \ \forall k=2,\dots,K  \\
    &   \bm{x}_{k} \leq \bm{u}_{k}\odot\bm{z}_{k}, \ \forall k=2,\dots,K  \\
    &   \bm{x}_{k} \geq \bm{W}_{k-1} \bm{x}_{k-1} + \bm{b}_{k-1}, \ \forall k=2,\dots,K  \\
    &   \bm{x}_{k} \geq 0, \ \forall k=2,\dots,K \\
    &   \bm{z}_{k} \in \{0,1\}^{n_k},  \ \forall k=2,\dots,K \\
    & \bm{x}_1 \in \mathcal{P},
    \end{array}
    \label{eq:multi_layer}
\end{equation}

\noindent where $\odot$ denotes the operator of the Hadamard product between two vectors. The upper and lower bounds of each neuron of layer $k$, namely $\bm{u}_k$ and $\bm{l}_k$, are calculated sequentially using LP. More specifically, $\bm{u}_k = \sup\limits_{ \bm{l}_{k-1} \leq \bm{x}_{k-1} \leq\bm{u}_{k-1}} \{ \bm{W}_{k-1}\bm{x}_{k-1} + \bm{b}_{k-1} \}$ and $\bm{l}_k = \inf\limits_{ \bm{l}_{k-1} \leq \bm{x}_{k-1} \leq \bm{u}_{k-1}} \{ \bm{W}_{k-1}\bm{x}_{k-1} + \bm{b}_{k-1} \}$ for all $k=3,\dots,K$. For $k=2$, $\bm{u}_2 = \sup\limits_{ \bm{x}_1 \in \mathcal{P}} \{ \bm{W}_{1}\bm{x}_{1} + \bm{b}_{1} \}$ and $\bm{l}_2 = \inf\limits_{ \bm{x}_1 \in \mathcal{P}} \{ \bm{W}_{1}\bm{x}_{1} + \bm{b}_{1} \}$.

Both in \eqref{eq:1_layer_mip} and in \eqref{eq:multi_layer} we come across MIP formulations. Despite their power in terms of expressivity, MIP formulations usually do not scale well as the size of the network increases. In what follows, we describe two sampling-based approaches that can solve the problem \eqref{eq:multi_layer} in higher dimensions.

\section{Sampling Approximation Algorithms}\label{sec:sampling_methods}

\subsection{Sampling Algorithm}\label{sec:sampling_method}

\par
\begin{figure}[ht]
      \centering
      \includegraphics[scale=0.1]{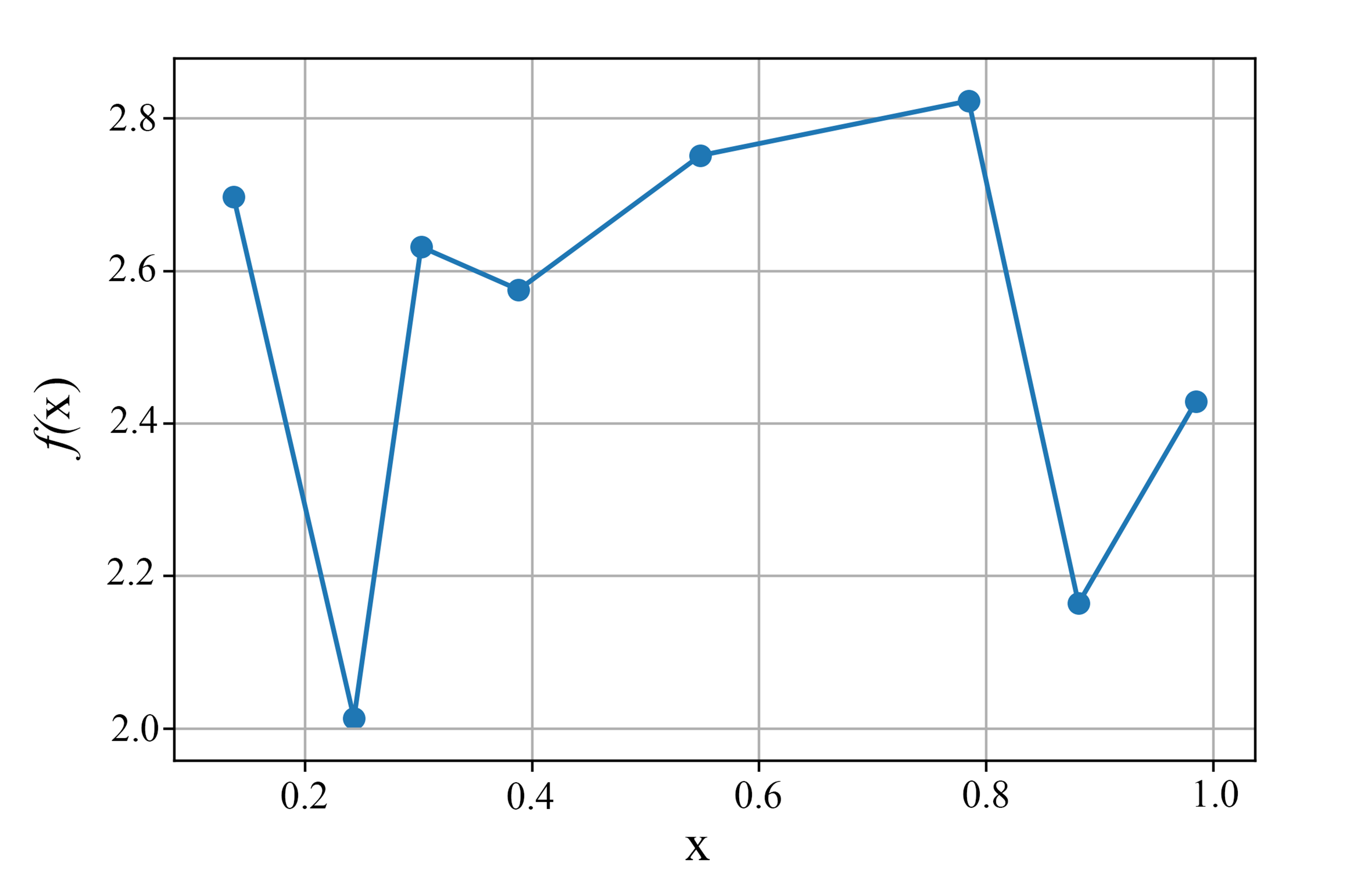}
      \caption{The output space of a ReLU neural network. The optimal solution is on a corner point.}
      \label{fig:1d}
\end{figure}

 We first propose a sampling-based algorithm that efficiently optimizes neural networks by reducing the MIP problem into multiple LPs via sampling. The algorithm takes advantage of two properties of ReLU neural networks. First, ReLU neural networks are piecewise linear functions and, as an immediate consequence, the optimal solution is always found at the extreme point (endpoint) of one of the output's hyperplanes. This is depicted in Figure \ref{fig:1d}, where we demonstrate a randomly initialized neural network in which the input $x$ is one dimensional and belongs to $\mathcal{P}=[0,1]$. 

\par Second, by using a ReLU neural network we can retrieve the hyperplane of the output space to which each input vector maps using formulation \eqref{eq:multi_layer}. Let $\bm{x}$ be the input vector sampled from $\mathcal{P}\subset \mathbb{R}^d$. By performing the feed-forward pass, we obtain which neurons are activated in the network. Therefore, we can set the corresponding binary variables of formulation \eqref{eq:multi_layer} to $1$ or $0$, depending on whether a neuron is activated or not. By setting the corresponding binary variables, the problem of finding the extreme point of the corresponding hyperplane reduces to an LP. 

\begin{figure}[ht]
    \centering
    \includegraphics[scale=0.1]{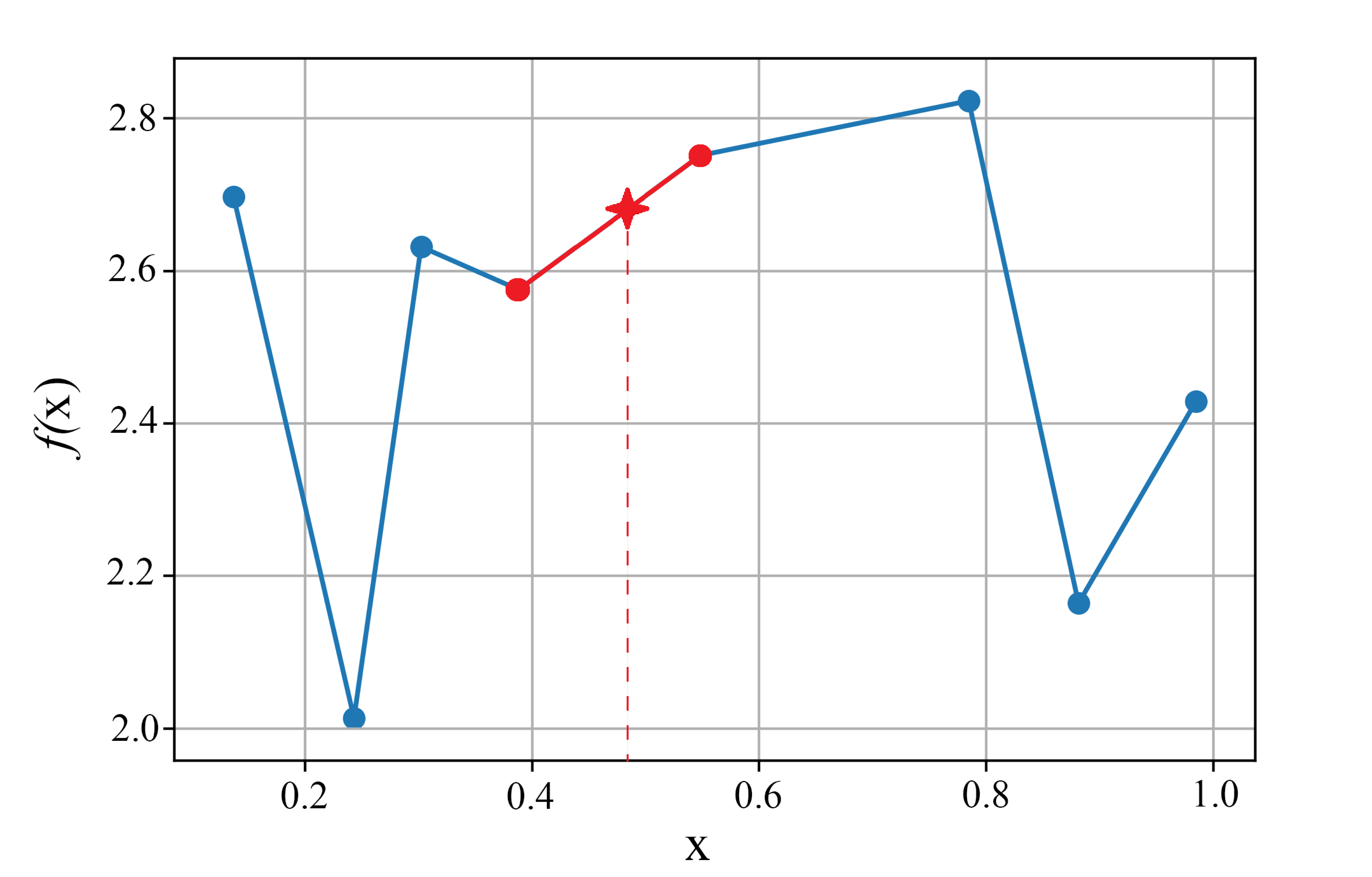}
    \caption{The red cross corresponds to the value of the network for the randomly selected input. The hyperplane (red line) to which the input vector maps is obtained via the forward pass and, by solving the LP over it, we obtain the optimal solution for this hyperplane (upper right red endpoint).}
    \label{fig:random_sample}
\end{figure}

This is depicted in Figure \ref{fig:random_sample}. The sampled input maps to the red cross in the output space of the network, and the corresponding hyperplane is the red line. By solving the LP over the hyperplane, we obtain that the optimal solution for this hyperplane is on its upper right endpoint. We can perform this procedure iteratively for multiple samples, to explore as many hyperplanes as possible and, retrieve an approximation of the global maximum.
\par More formally, if we denote as $f(\bm{x};\bm{z})$ the neural network output with fixed activated neurons, i.e., fixed binary variables $\bm{z}$, our approximation of the global maximum after $N$ samples is given by the following expression:
\begin{align}
    \max\{ \max\limits_{\bm{x} \in \mathcal{P}}f(\bm{x};\bm{z}^1),\dots,\max\limits_{\bm{x} \in \mathcal{P}}f(\bm{x};\bm{z}^N) \}, \label{eq:compact}
\end{align}

\noindent where $\bm{z}^i, \ i \in \{1,\dots,N\}$, corresponds to the activated neurons when the $i$-th sampled input vector is evaluated. Equation \eqref{eq:compact} can be rewritten as the following MIP.
\begin{equation}
    \begin{array}{ll@{}ll}
    \max\limits_{\bm{x},\bm{z}} & f(\bm{x};\bm{z}) \\
    \text{subject to}   
    &  \bm{z} \in \{\bm{z}^1,\dots,\bm{z}^N\} \\
    &  \bm{x} \in \mathcal{P}. \\
    \end{array}
    \label{eq:total_mip}
\end{equation}

\par Problem \eqref{eq:total_mip} can be decoupled into $N$ distinct LPs since $\bm{x}$ is a continuous variable, $f$ is a ReLU neural network and $\bm{z}$ can take only $N$ discrete values. Based on this observation, we come up with the algorithm described in pseudocode \ref{alg:sampling}. Given the number of samples $N$, the ReLU neural network $f$, and the domain of the input $\mathcal{P}$, the algorithm for each sampled input vector $\bm{x}^i$ performs the feed-forward pass, retrieves $\bm{z}^i$, fixes it in \eqref{eq:multi_layer}, and then solves the LP. The algorithm returns as the optimal solution the extreme point $\bm{x}_{max}$ that achieves the highest objective value after $N$ iterations. Since for each sample, we solve an LP, this method takes advantage of the scalability of the existing LP solvers and scales for larger neural networks than MIP-based approaches do.

\begin{algorithm}

\caption{Sampling Algorithm}
\begin{algorithmic}[1]
\STATE \textbf{INPUT}: $N, f, \mathcal{P}$
\STATE $i \leftarrow 0$, $max \leftarrow -INF$,  $\bm{x}_{max} \leftarrow None$
\WHILE{$i<N$} 
\STATE Sample $\bm{x}^i$ from $\mathcal{P}$.
\STATE $\bm{z}^i \leftarrow get\_activations(f,\bm{x}^i)$ // \begin{footnotesize}\texttt{ performs the feed-forward pass for input $\bm{x}^i$ and retrieves the corresponding activated neurons.}\end{footnotesize}
\STATE $(max\_lp,\bm{x}^*) \leftarrow solve\_lp(f,\bm{z}^i) $ // \begin{footnotesize}\texttt{solves \eqref{eq:multi_layer} with fixed $\bm{z}^i$, returns optimal objective and solution.}\end{footnotesize}
\IF{ $max\_lp > max$} 
\STATE $max \leftarrow max\_lp$
\STATE $\bm{x}_{max} \leftarrow \bm{x}^*$
\ENDIF
\STATE $i \leftarrow i+1$
\ENDWHILE
\RETURN ($\bm{x}_{max}, max$)
\end{algorithmic}
\label{alg:sampling}
\end{algorithm}

\subsection{Enhanced Sampling Algorithm}\label{sec:local}

In the algorithm discussed above, after solving the LP, the algorithm restarts by sampling a new input vector $\bm{x}$ from the feasible region $\mathcal{P}$. Nevertheless, the algorithm does not explore the landscape around the solution to which it gives rise at each iteration. For instance, in the example depicted in Figure \ref{fig:random_sample}, after finding the local optimal solution on the upper right corner point of the red hyperplane, instead of re-sampling, it would be beneficial to first follow the positive slope at this corner point, move to the next hyperplane and obtain an improved solution. In what follows, we describe how to enhance the sampling algorithm using local search.

Regarding the inputs of the enhanced algorithm, they remain the same with the addition of a $\textit{gap}$ parameter that controls the trade-off between exploration and exploitation. The enhanced algorithm, as before, samples an input vector $\bm{x}^i$ from $\mathcal{P}$, performs the feed-forward pass, retrieves the vector $\bm{z}^i$, fixes it in \eqref{eq:multi_layer}, and solves the corresponding LP. Then, the algorithm calculates the absolute relative difference between the objective of the solution of the LP and the current highest objective. If the absolute relative difference is less than the allowed $\textit{gap}$, i.e., the retrieved solution is not far from the current highest objective, the algorithm performs a local search, where the algorithm searches for neighboring hyperplanes with corner points that can lead to higher objectives.
 
Local search is based on the following observation. Each solution of the LP often belongs to more than one hyperplane, as it is a corner point in a piecewise linear function. In ReLU neural networks, corner points correspond to neurons that output a value before the ReLU activation equal to $0$ (in practice, $\epsilon$-close to $0$). Therefore, by having as unknowns in \eqref{eq:multi_layer} only the binary values that correspond to the neurons that output a zero value before the activation, and by fixing the others as before, we can solve a MIP that is small in terms of binary variables and constraints and which, at the same time, optimizes over the neighboring hyperplanes, that can lead to a possibly improved solution. Local search is repeated until the best, in terms of objective value, neighboring hyperplane is found.

Regarding the trade-off between exploration and exploitation that the $\textit{gap}$ parameter introduces, for low values of $\textit{gap}$ we begin the local search with corner points that have objective values close to the current highest, while for high values we explore corner points that have a lower objective. In the former case, we solve smaller MIPs and, consequently, the method is more efficient computationally. Nevertheless, this also implies that the method explores fewer hyperplanes. In the latter case, we solve larger MIPs, which is less efficient computationally, but at the same time, it allows the algorithm to explore more neighboring hyperplanes. For the computations, we picked $\textit{gap}=\frac{2}{3}$. A description of the algorithm can be found in pseudocode \ref{alg:sampling_local}.

\begin{algorithm}

\caption{Enhanced Sampling Algorithm}
\begin{algorithmic}[1]
\STATE \textbf{INPUT}: $N, f, \mathcal{P}, \textit{gap}$
\STATE $i \leftarrow 0$, $max \leftarrow -INF$, $\bm{x}_{max} \leftarrow None$
\WHILE{$i<N$} 
\STATE Sample $\bm{x}^i$ from $\mathcal{P}$.
\STATE $\bm{z}^i \leftarrow get\_activations(f,\bm{x}^i)$
\STATE $(max\_lp,\bm{x}^*) \leftarrow solve\_lp(f,\bm{z}^i) $ 
\IF{ $max\_lp > max$} 
    \STATE $max \leftarrow max\_lp$
    \STATE $\bm{x}_{max} \leftarrow \bm{x}^*$
\ENDIF
\IF{ $|(max-max\_lp)/max| < \textit{gap}$} 
    \WHILE {\textit{True}}
        \STATE $(max\_ip,\bm{x}^*) \leftarrow solve\_mip(f,\bm{x}^*) $ // \begin{footnotesize}\texttt{solves \eqref{eq:multi_layer} for $\bm{x}^*$ with unknowns only the binary variables of neurons with $0$ value before the ReLU.}\end{footnotesize}
        \IF{ $max\_ip > max\_lp$} 
            \STATE $max\_lp \leftarrow max\_ip$
            \IF{ $max\_lp > max$} 
                \STATE $max \leftarrow max\_lp$
                \STATE $\bm{x}_{max} \leftarrow \bm{x}^*$
            \ENDIF
        \ELSE
            \STATE \textbf{break}
        \ENDIF    
    \ENDWHILE
\ENDIF
\STATE $i \leftarrow i+1$
\ENDWHILE
\RETURN ($\bm{x}_{max}, max$)
\end{algorithmic}
\label{alg:sampling_local}
\end{algorithm}

\subsection{Provable Guarantees}\label{sec:proovable_guarantees}

\par In this section, we prove the convergence and the sample complexity of the proposed methods and we show that the enhanced sampling algorithm finds a local optimum from the first iteration. Since we are optimizing over ReLU neural networks over bounded domains, there is a finite number of partitions of the input space to hyperplanes (\cite{montufar2014number, serra2018bounding}). Every time we uniformly sample a point from the domain, we fall into one of those hyperplanes, and after optimizing, we obtain the corresponding maximum value for the specific hyperplane. We denote as $X_i$, for $i=1,\dots, N$, the random variable indicating the maximum value obtained after optimizing over the hyperplane to which sample $i$ corresponds. The family of $X_i$ are i.i.d. random variables.

\begin{lemma}(\textit{von Mises' condition}, \cite{von1936distribution}) Let $F$ be a CDF and $M\in \mathbb{R}\cup\{+\infty\}$ its right endpoint. Suppose that $F''(x)$ exists and $F'(x)$ is positive for all $x$ in a neighborhood of $M$. If 
\begin{align}
    \lim_{t \to M } \Big( \frac{1-F}{F'}\Big)'(t)  = \gamma, \textit{ or equivalently,} \lim_{t \to M } \Big( \frac{(1-F(t))F''(t)}{(F'(t))^2}\Big)  = -\gamma-1, 
\end{align}
then $F$ satisfies von Mises' condition.
\end{lemma}
\begin{theorem}\label{thm:proof}
\textit{Algorithm \ref{alg:sampling} has the following properties: }
\begin{enumerate}

    \item\textit{ (\textbf{Convergence}) Algorithm \ref{alg:sampling} converges to the optimal solution of \eqref{eq:multi_layer} almost surely.}
     \item \textit{(\textbf{Sample Complexity}) Assume that the CDF $F(\cdot)$ of $X_i$'s satisfies von Mises' condition. Let $\epsilon >0$ be appropriately small and be $N$ the number of samples so that:}
    \begin{align}
        N \geq \frac{\ln(\frac{1}{\delta})}{\epsilon^{-\frac{1}{\gamma}}}.
    \end{align}
    \textit{Then, Algorithm \ref{alg:sampling} using $N$ samples is $\epsilon$-close to the true maximum $M$ with probability of at least $1 - \delta$.}
\end{enumerate}

\end{theorem}

The theorem also applies to Algorithm \ref{alg:sampling_local} and the complete proof is in Section~\ref{sec:proof} of the Appendix.
\par To calculate the exact lower bound on the sample complexity, we need to understand two things. First, we need to verify that the CDF satisfies von Mises'  condition and, second, we need to have knowledge of the parameter $\gamma$, also called the \textit{extreme value index}. In Section~\ref{sec:vonmises} of the Appendix, we show computationally how to verify von Mises' condition and estimate the extreme value index $\gamma$.

\begin{remark}
  \textit{Algorithm \ref{alg:sampling_local} reaches a local optimum from the first iteration.}
\end{remark}
According to Algorithm \ref{alg:sampling_local}, to sample a new point, the break condition in line 21 should be satisfied. This only happens when local search cannot find a solution that increases the current objective in any neighboring hyperplane of the current solution. Therefore, the break condition is satisfied only when the current solution is a local optimum. That means that the enhanced sampling algorithm reaches a local optimum from the first iteration.


\section{Computational Experiments}\label{sec:experiments}

Computational experiments were performed using Gurobi v 9.5 (\cite{gurobi}, MIT License) and Julia 1.5.2 (\cite{Julia-2017}, MIT License). All experiments were performed on an internal cluster with a 2.20GHz Intel(R) Xeon(R) Gold 5120 CPU and 256 GB memory.

\subsection{Comparison with other methods}\label{sec:comparison}
To measure the performance of the sampling algorithms, we evaluate them over multiple synthetic feed-forward ReLU neural networks against other state-of-the-art methods. More specifically, we report the maximum objective value achieved by the sampling algorithm, the enhanced sampling algorithm, the improved MIP formulation by \cite{fischetti2018deep}, the Big-M + cuts formulation by \cite{anderson2020strong} and the $N=\{2,4\}$ equal-size partition by \cite{tsay2021partition}.

We evaluate all methods on one- and two-layer feed-forward neural networks with $100$ neurons per layer. The weights of the networks are randomly initialized according to the uniform Xavier initialization method (\cite{glorot2010understanding}).  The input dimension $d$ varies from $5$ to $1,000$ and the domain of $\mathcal{P}$, from which we sample uniformly at random, is $[0,1]^d$. To make a fair comparison, all methods were given a one-hour time window. The results are shown in Tables \ref{tab:one_layer} and \ref{tab:two_layer}.

\begin{table}[ht]
    \centering
        \caption{Maximum objective value achieved by each method for the 1-layer case.}
    \label{tab:one_layer}
\begin{tabular}{cccccccc}
\toprule 
     \textbf{$d$}  & Sampling & Enhanced Sampling& MIP & Big-M + cuts & 2 Partitions & 4 Partitions \\ 

    \midrule

    5 &   \textbf{-0.0721} &      \textbf{-0.0721} & \textbf{-0.0721} &       \textbf{-0.0721} &     \textbf{-0.0721} &     -\textbf{0.0721} \\
   10 &    0.3738 &        \textbf{0.3813} &  \textbf{0.3813} &        \textbf{0.3813} &      \textbf{0.3813} &      \textbf{0.3813} \\
   20 &    1.0584 &       \textbf{1.1021} &  \textbf{1.1021} &        1.0959 &      \textbf{1.1021} &      \textbf{1.1021 }\\
   50 &    2.1369 &       \textbf{2.2086} & \textbf{2.2086} &        2.2053 &      2.2000 &      2.1962 \\
  100 &    3.3274 &        \textbf{3.6144} &  \textbf{3.6144} &        \textbf{3.6144} &      3.5279 &      3.5496 \\
   200 &    3.1809 &        \textbf{3.2863} &  3.2669 &        3.2149 &      3.0165 &      3.2623 \\
   500 &    3.7571 &        \textbf{4.0007} &  3.9544 &        3.9907 &      3.7820 &      3.6084 \\
  750 &    4.5920 &        \textbf{4.9097} &  \textbf{4.9097} &        4.8949 &      4.8041 &      4.6281 \\
 1000 &   4.4134  &  \textbf{4.5509}      &  4.5428  & \textbf{4.5509}        & 4.5381     & 4.4404      \\ 
\bottomrule
\end{tabular} 
\end{table}

\begin{table}[ht]
    \centering
        \caption{Maximum objective value achieved by each method for the 2-layer case. \texttt{NA} means that the method did not manage to produce a solution during the given time window.}
    \label{tab:two_layer}
\begin{tabular}{cccccccc}
\toprule 
     \textbf{$d$}  & Sampling & Enhanced Sampling& MIP & Big-M + cuts & 2 Partitions & 4 Partitions \\ 

    \midrule

    5 &    \textbf{0.0627} &        \textbf{0.0627} &  \textbf{0.0627} &        0.0307 &      0.0242 &      0.0304 \\
   10 &    0.3361 &       \textbf{0.3680} & \textbf{0.3680} &        0.2915 &         \texttt{NA} &        \texttt{NA} \\
   20 &   0.0090  &  \textbf{0.0869}     &  0.0714 &   -0.0225    &    -0.0299  &   0.0204    \\
   50 &    0.9253 &        \textbf{1.1839} &  1.0803 &        0.8835 &      0.7931 &      0.8641 \\
   100 &    1.4630 &        \textbf{1.8462} &  1.7687 &        1.7346 &        \texttt{NA} &         \texttt{NA} \\
   200 &    1.8437 &        \textbf{2.5732} &  2.4605 &        2.1301 &      1.5630 &      2.0957 \\
   500 &    3.8582 &        \textbf{5.0336} &  4.7706 &         \texttt{NA} &      2.2235 &      2.8019 \\
   750 &    4.6663 &        \textbf{5.3842} &  5.1609 &        4.7887 &      3.1149 &     -0.0262 \\
  1000 &    6.1788 &        \textbf{7.4953} &  5.9709 &        6.0861 &      5.1820 &      3.2248 \\
\bottomrule
\end{tabular} 
\end{table}

In Table \ref{tab:one_layer}, we observe that for one-layer networks and low input dimensions ($5-10$) all MIP-based methods, as well as the enhanced sampling algorithm, reach the same maximum objective value. As the input dimension increases ($d \geq 20$), we observe that the enhanced sampling algorithm performs at least as well as the rest of the methods, with the improved MIP and the Big-M + cuts methods performing a bit worse in some cases, but similarly well overall. For the two-layer case, we observe that as the input dimension increases, the MIP-based methods begin to perform significantly worse compared to the enhanced sampling algorithm. More specifically, as shown in Table \ref{tab:two_layer}, for input dimensions $d=5,10$ the enhanced sampling algorithm and the improved MIP method perform equally well, while for input dimension $d$ higher than $10$ the enhanced sampling algorithm performs the best, by reaching approximately $4.3\%-25.0\%$ higher objective value than the second-best performing method, depending on the case. It is worth noting that there are cases, where some MIP-based methods did not manage to produce a solution during the given time window.

Overall, the enhanced sampling algorithm performs at least as well as the best-performing MIP-based methods for the one-layer case, while for the two-layer case, it performs consistently better than any other MIP-based method for input dimension $d$ higher than $10$. Next, we show how the sampling algorithms can be used to warm-start the MIP-based formulations and improve their performance.

\subsection{Sampling as a Warm Start}\label{sec:warm}

\par  As shown in Table \ref{tab:two_layer}, there are cases where the MIP-based formulations cannot produce a solution during the given time window. This usually happens because the solver is slow in finding an initial feasible solution. One solution is to use one of the sampling algorithms to produce an initial feasible solutions as a warm-start for the solver. A benefit of this approach is that it is computationally efficient, since each iteration of the sampling algorithms corresponds to solving a single LP (and in the enhanced sampling method, maybe a few, small in size, MIPs).  

To evaluate this approach, we conduct a second round of experiments by producing warm-starts using the initial sampling algorithm on the same networks. We select the initial sampling algorithm for producing warm-starts, instead of the enhanced sampling algorithm, as it is extremely fast in practice. As before, all methods were given a one-hour time window. To produce the warm-start solution for each network, we run the sampling algorithm for $1,000$ iterations. The results are shown in Tables \ref{tab:one_layer_warm} and \ref{tab:two_layer_warm}.

\begin{table}[ht]
    \centering
        \caption{Maximum objective value achieved by each method for the 1-layer case using warm-starts.}
    \label{tab:one_layer_warm}
\begin{tabular}{cccccccc}
\toprule 
     \textbf{$d$}  & Sampling & Enhanced Sampling& MIP & Big-M + cuts & 2 Partitions & 4 Partitions \\ 

    \midrule

  5 &   \textbf{-0.0721} &       \textbf{-0.0721} & \textbf{-0.0721} &       \textbf{-0.0721} &     \textbf{-0.0721} &     \textbf{-0.0721} \\
    10 &    0.3744 &        \textbf{0.3813} &  \textbf{0.3813} &        \textbf{0.3813} &      \textbf{0.3813} &      \textbf{0.3813} \\
   20 &    1.0748 &        \textbf{1.1021} &  \textbf{1.1021} &       \textbf{1.1021} &     \textbf{1.1021} &      \textbf{1.1021} \\
   50 &    2.1362 &        \textbf{2.2086} &  2.2053 &        \textbf{2.2086} &      2.1895 &      2.2008 \\
   100 &    3.3927 &        \textbf{3.6144} & \textbf{ 3.6144} &        3.6142 &      3.5794 &      3.6138 \\
   200 &    3.1935 &        \textbf{3.2863} &  3.2807 &        3.2122 &      3.2146 &      3.2693 \\
  500 &    3.7362 &        \textbf{4.0007} &  3.9940 &        3.9947 &      3.7259 &      3.8075 \\
   750 &    4.5813 &        \textbf{4.9097} &  \textbf{4.9097} &        4.8777 &      4.8088 &      4.6451 \\
   1000 &   4.4134  &  \textbf{4.5509}      &  \textbf{4.5509}  & 4.5438       & 4.4697     & 4.4804      \\
  
\bottomrule
\end{tabular} 
\end{table}

\begin{table}[ht]
    \centering
        \caption{Maximum objective value achieved by each method for the 2-layer case using warm-starts.}
    \label{tab:two_layer_warm}
\begin{tabular}{cccccccc}
\toprule 
     \textbf{$d$}  & Sampling & Enhanced Sampling& MIP & Big-M + cuts & 2 Partitions & 4 Partitions \\ 

    \midrule

  5 &    \textbf{0.0627}&        \textbf{0.0627} &  \textbf{0.0627} &        0.0371 &      0.0373 &      0.0373 \\
  10 &    0.3376 &        \textbf{0.3680} &  0.3365 &        0.3313 &      0.3367 &      0.3319 \\
  20 & 0.0017  &    \textbf{0.0869}    & 0.0727  & -0.0467 & 0.0355 &  0.0195 \\
   50 &    0.8973 &        \textbf{1.1839} &  1.1168 &        1.0480 &      1.0673 &      1.0383 \\
 100 &    1.4688 &        \textbf{1.8462} &  1.7938 &        1.6022 &      1.3504 &      1.3504 \\
  200 &    1.8446 &        \textbf{2.5569} &  2.4134 &        1.8683 &      1.8233 &      1.8233 \\
  500 &    4.0620 &        \textbf{5.0024} &  4.8766 &        3.7529 &      3.7529 &      3.7529 \\
 750 &    4.6945 &        \textbf{5.4139} &  5.0903 &        4.7887 &      4.5359 &      4.5359 \\
 1000 &    5.8908 &        \textbf{7.5187} &  6.9233 &        6.0861 &      5.5857 &      5.5857 \\
\bottomrule
\end{tabular} 
\end{table}

In Table \ref{tab:one_layer_warm}, we observe that for one-layer networks warm-starts slightly increase the objective values reached (increase between $0.16\%$ and $1.29\%$) from MIP-based methods. This happens since for the one-layer case, the network is small in terms of parameters and the MIP-based methods can produce strong solutions without warm-starts, as shown in Table \ref{tab:one_layer}. For the two-layer case, as shown in Table \ref{tab:two_layer_warm}, we observe that warm-starts can significantly improve the performance of MIP-based methods. More specifically, warm-starts for the specific experiment improved on average the improved MIP method by $4.07\%$, the Big-M + cuts method by $22.57\%$, and the $N=\{2,4\}$ equal-size partition method by $63.75\%$ and $2,511\%$ accordingly. More importantly, warm-starts were able to tackle the problem of not producing a solution during the given time window, as illustrated in Table \ref{tab:two_layer}. 

Overall, even though improvements like the ones reported for the $N=\{2,4\}$ equal-size partition may seem ambitious and case-dependent, using the sampling algorithms to produce initial feasible solutions as a warm-start for the solver seems to be beneficial, especially for larger, in terms of parameters, neural networks.

\section{Conclusions}\label{sec:conclusions}

This paper introduces efficient, sampling-based algorithms for optimizing objective functions that come from trained ReLU neural networks. We first propose an iterative algorithm that takes advantage of the piecewise linear structure of ReLU neural networks and reduces the initial MIP into multiple easy-to-solve LPs via sampling. We then improve upon it by introducing an enhanced method that considers at each iteration the neighborhood of the solution of the LP and, reduces the initial MIP problem into a few, smaller, and easier-to-solve MIPs. Even though for the ease of exposition, the computations were conducted over one type of neural network with one and two layers, experimental evidence shows that the proposed algorithms outperform current, MIP-based, state-of-the-art methods and can be used efficiently to warm-start MIP-based methods, especially for larger, in terms of parameters, neural networks.

\textbf{Broader Impact Statement}: Our algorithms can simplify and improve the process of optimizing objective functions from trained ReLU neural networks, which can potentially improve tasks such as the verification of the robustness of a neural network output, the compression of neural networks, and the retrieval of adversarial examples among others. Since it is a fundamental optimization methodology, we do not foresee a negative impact on the society implied by the algorithms directly.

\bibliographystyle{plainnat}
\bibliography{references}

\newpage
\appendix

\section{Proof of Theorem \ref{thm:proof}}\label{sec:proof}

According to Theorem \ref{thm:proof}: 
\begin{enumerate}

    \item (\textbf{Convergence}) Algorithm \ref{alg:sampling} converges to the optimal solution of \eqref{eq:multi_layer} almost surely.
     \item (\textbf{Sample Complexity}) Assume that the CDF $F(\cdot)$ of $X_i$'s satisfies von Mises' condition. Let $\epsilon >0$ be appropriately small and be $N$ the number of samples so that:
    \begin{align}
        N \geq \frac{\ln(\frac{1}{\delta})}{\epsilon^{-\frac{1}{\gamma}}}.
    \end{align}
    Then, Algorithm 1 using $N$ samples is $\epsilon$-close to the true maximum $M$ with probability of at least $1 - \delta$.
\end{enumerate}

\begin{proof}

1. (\textbf{\textit{Convergence}}) Let $N>0$. We define: 
\begin{align}
    M_N := \max \{X_1,\dots,X_N\}.
\end{align}
We have that 
\begin{align}
    \mathbb{E}[M_{N+1}-M_N|M_N] = \mathbb{E}[M_{N+1}|M_N] -M_N \geq 0,
\end{align}
\noindent which means that $M_N$ is a submartingale. Since we are maximizing over a bounded domain, we have that $|M_N|<\infty$ for every $N$. Therefore,  $\sup_N \mathbb{E}[|M_N|] < \infty$. As a result, according to Doob's martingale convergence theorem, $M_N$ converges to the maximum value, denoted by $M$, almost surely. $\square$

2. (\textbf{\textit{Sample Complexity}}) Let $\epsilon >0$ be the gap parameter. We have that:
\begin{align*}
    \mathbb{P}[ M-M_N \geq \epsilon  ] &= \mathbb{P}[ M_N \leq M - \epsilon  ] = \mathbb{P}[ \max\{X_1,\dots,X_N\} \leq M - \epsilon ] =\\
    &= \mathbb{P}[ X_1 \leq M - \epsilon,\dots,X_N \leq M - \epsilon ] \myeq \mathbb{P}[ X_1 \leq M - \epsilon ]^N. \numberthis \label{eq:bound1}
\end{align*}
From \eqref{eq:bound1}, we have that
\begin{align}
    \mathbb{P}[ M-M_N \geq \epsilon  ] &= \mathbb{P}[ X_1 \leq M - \epsilon ]^N.  \label{eq:bound2}
\end{align}

\noindent Since $F(\cdot)$ satisfies von Mises' condition, for $\epsilon >0$ appropriately small, we obtain the following (see \cite{de2006extreme}):
\begin{align}
    \mathbb{P}[X_1 \leq M-\epsilon]^N = (F(M-\epsilon))^N \approx (1-\epsilon^{-\frac{1}{\gamma}})^N. \label{eq:bound3}
\end{align}

\noindent By using the inequality $e^{-x} \geq 1-x$, for $1 \geq x\geq0$, \eqref{eq:bound3} becomes
\begin{align*}
  \mathbb{P}[ X_1 \leq M - \epsilon ]^N &\leq \exp( -\epsilon^{-\frac{1}{\gamma}})^N , \label{eq:bound4} \numberthis
\end{align*}
\noindent We can set the right-hand side of the bound to be less than or equal to $\delta$ and then solve for $N$. Specifically,
\begin{align*}
    \exp( -\epsilon^{-\frac{1}{\gamma}})^N \leq \delta &\implies N(-\epsilon^{-\frac{1}{\gamma}}) \leq \ln (\delta)\implies N \geq \frac{\ln(\frac{1}{\delta})}{\epsilon^{-\frac{1}{\gamma}}}. \numberthis
\end{align*} 

\end{proof}

\section{Verifying  von Mises' condition}\label{sec:vonmises}

In what follows, we aim to illustrate how to verify von Mises' condition and estimate the extreme value index $\gamma$ computationally. To perform that, we run the sampling algorithm for a finite number of iterations to retrieve a set of $X_i$ instances. To obtain $X_i$'s for multiple neural networks, we conduct a series of experiments on one- and two-layer feed-forward ReLU neural networks with $100$ neurons at each layer. The weights of the networks are randomly initialized according to the uniform Xavier initialization method.  The input dimension $d$ varies from $5$ to $1,000$ and the domain of the input $\mathcal{P}$ is $[0,1]^d$. We run the sampling algorithm for $30$ minutes.

\par To verify computationally that the CDF of the family of $X_i$'s verifies von Mises' condition, we perform probability density fitting. The software that we use for this task is \texttt{distfit} (\cite{erdogant2019distfit}, MIT License). We first fit several well-known probability distributions on $X_i$'s, such as exponential, Pareto, Weibull, Student's t, F, generalized extreme value, and beta, and then, we select the distribution with the lowest residual sum of squares (RSS), while we present the top $3$ fitted distributions for each ReLU neural network. 

\begin{figure}[htp]
\centering

\includegraphics[width=.45\textwidth,height=0.2\textheight,keepaspectratio]{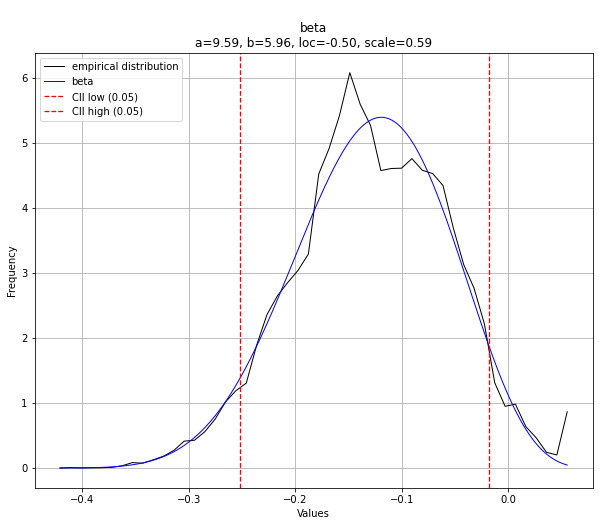}\quad
\includegraphics[width=.45\textwidth,height=0.2\textheight,keepaspectratio]{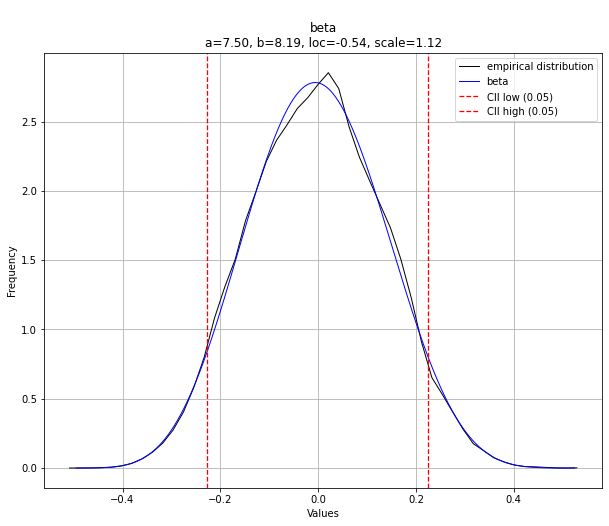}

\medskip

\includegraphics[width=.45\textwidth,height=0.2\textheight,keepaspectratio]{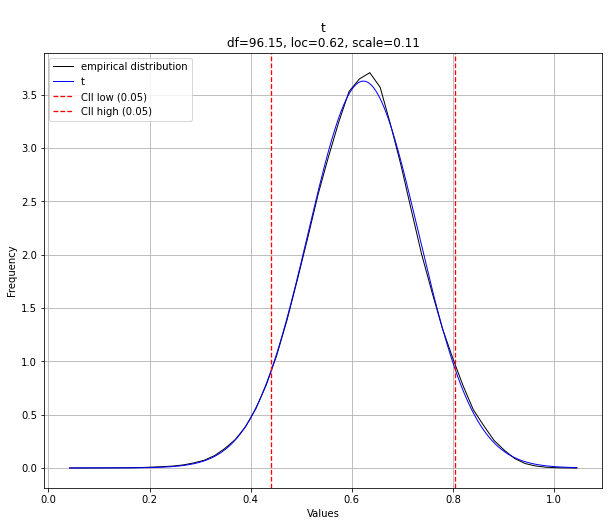}\quad
\includegraphics[width=.45\textwidth,height=0.2\textheight,keepaspectratio]{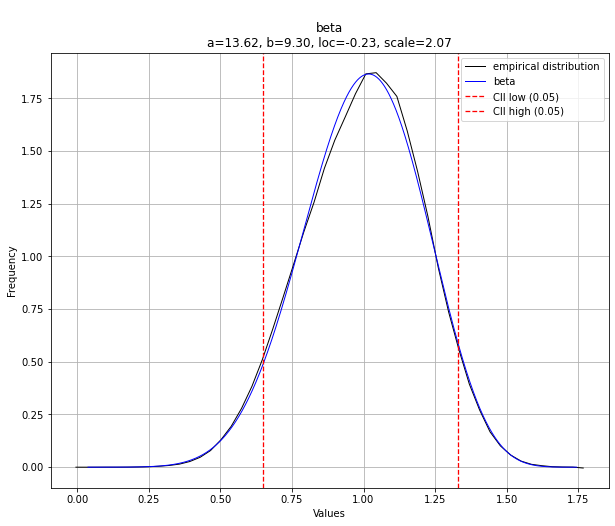}

\medskip

\includegraphics[width=.45\textwidth,height=0.2\textheight,keepaspectratio]{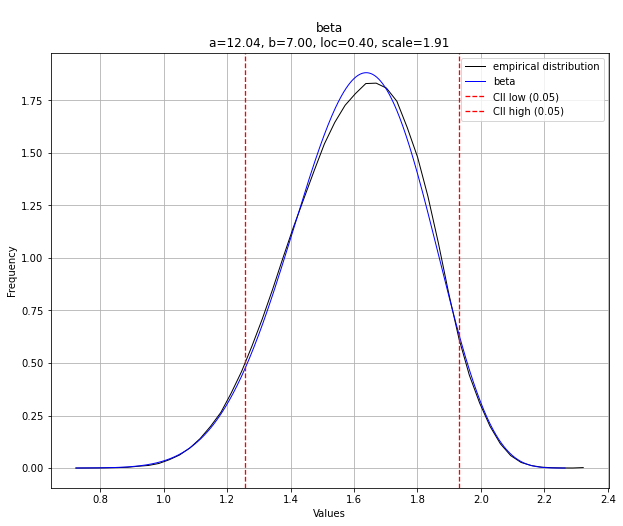}\quad
\includegraphics[width=.45\textwidth,height=0.2\textheight,keepaspectratio]{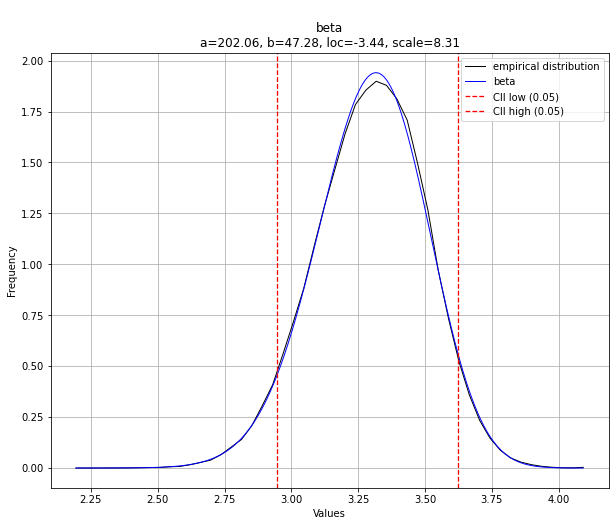}

\medskip

\includegraphics[width=.45\textwidth,height=0.2\textheight,keepaspectratio]{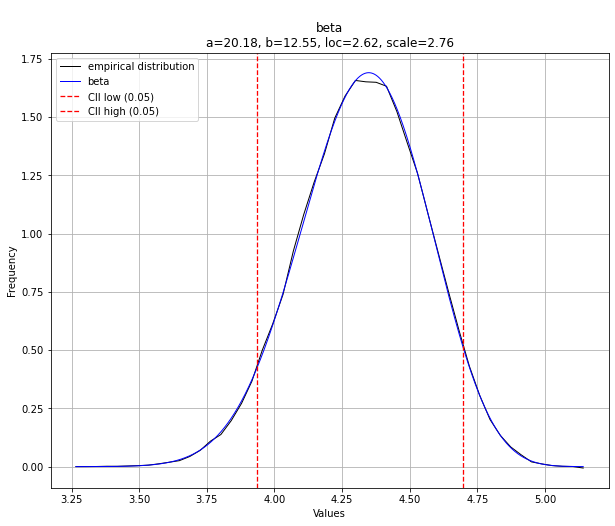}\quad
\includegraphics[width=.45\textwidth,height=0.2\textheight,keepaspectratio]{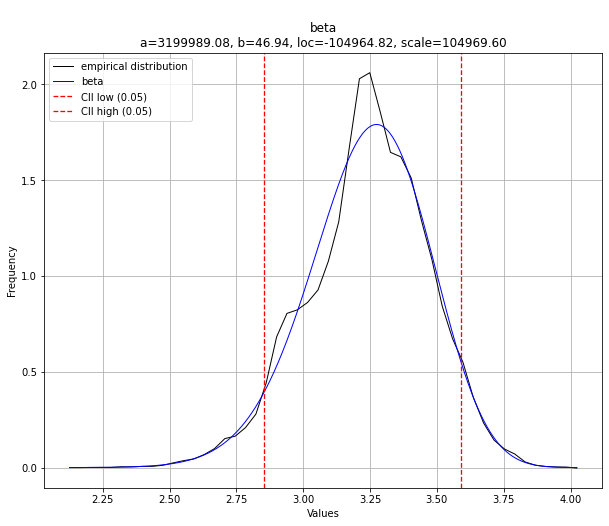}

\caption{Best-fitted probability distributions to the actual data for the one-layer case. (From left to right) First row: $d=5$ and $d=10$; Second row: $d=20$ and $d=50$; Third row: $d=100$ and $d=200$; Fourth row: $d=500$ and $d=1,000$.}
\label{fig:estimated_dist}
\end{figure}

\begin{figure}[htp]
\centering

\includegraphics[width=.45\textwidth,height=0.2\textheight,keepaspectratio]{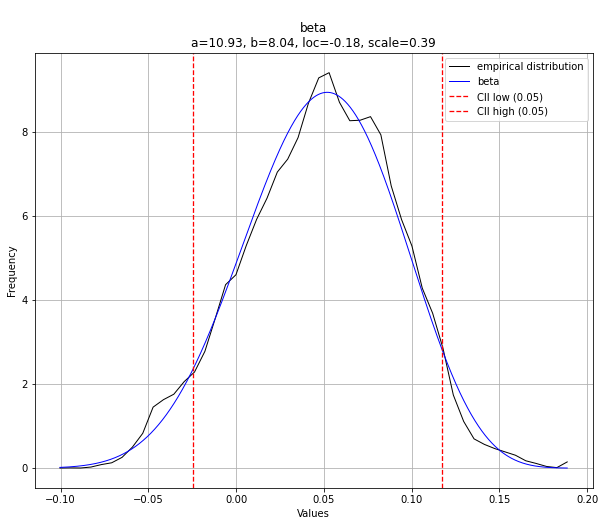}\quad
\includegraphics[width=.45\textwidth,height=0.2\textheight,keepaspectratio]{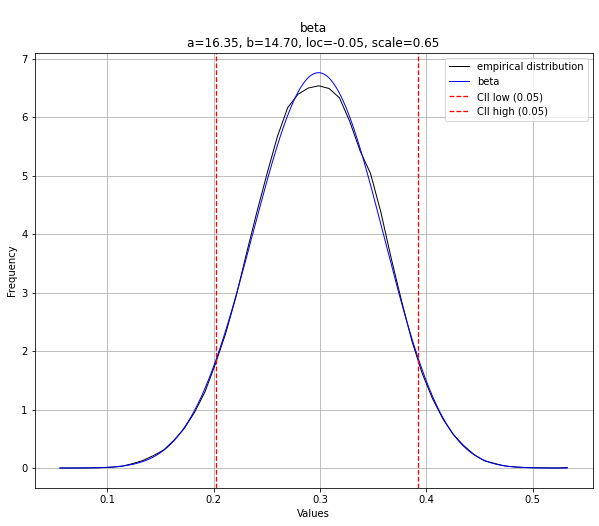}

\medskip

\includegraphics[width=.45\textwidth,height=0.2\textheight,keepaspectratio]{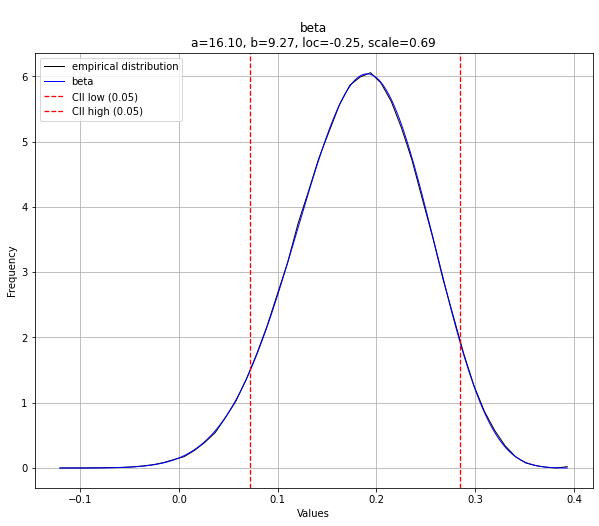}\quad
\includegraphics[width=.45\textwidth,height=0.2\textheight,keepaspectratio]{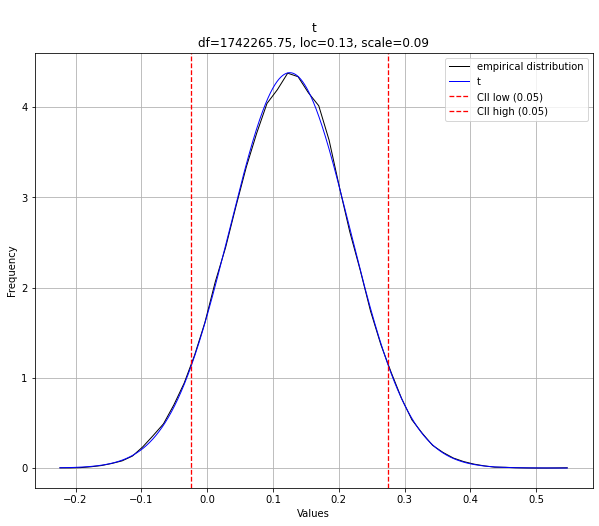}

\medskip

\includegraphics[width=.45\textwidth,height=0.2\textheight,keepaspectratio]{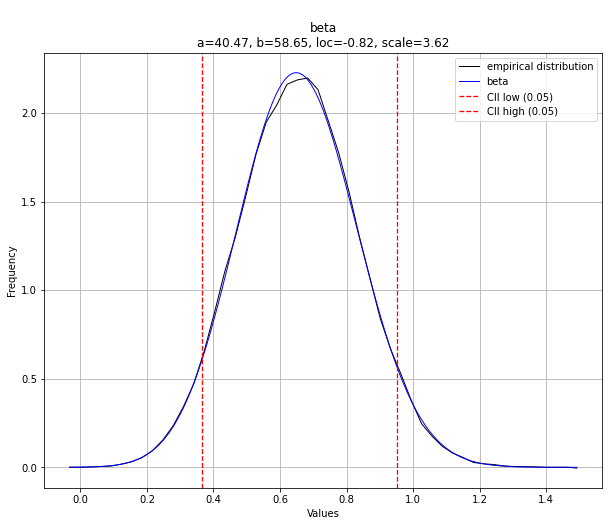}\quad
\includegraphics[width=.45\textwidth,height=0.2\textheight,keepaspectratio]{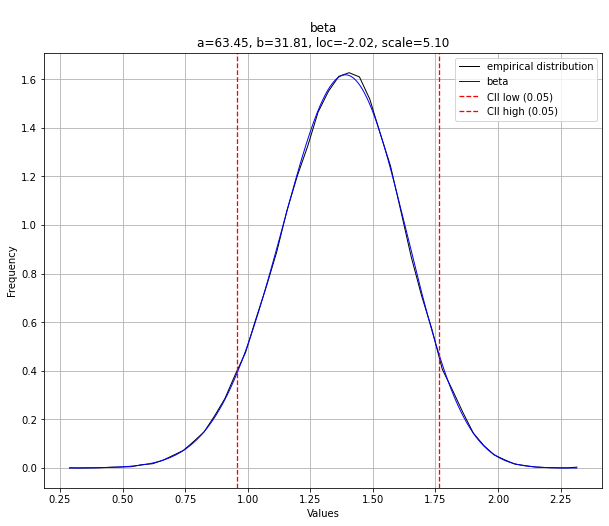}

\medskip

\includegraphics[width=.45\textwidth,height=0.2\textheight,keepaspectratio]{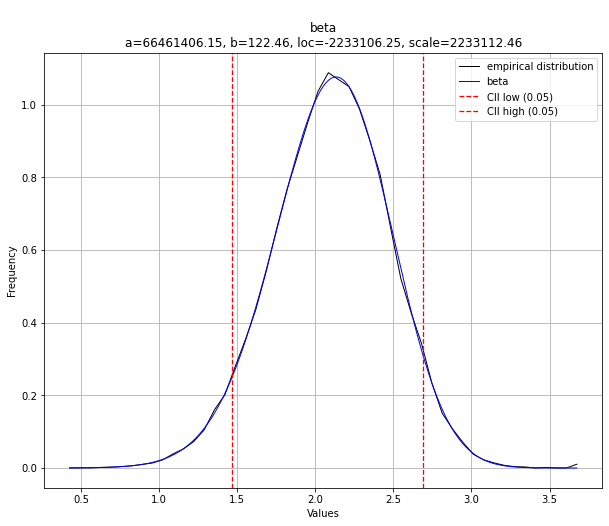}\quad
\includegraphics[width=.45\textwidth,height=0.2\textheight,keepaspectratio]{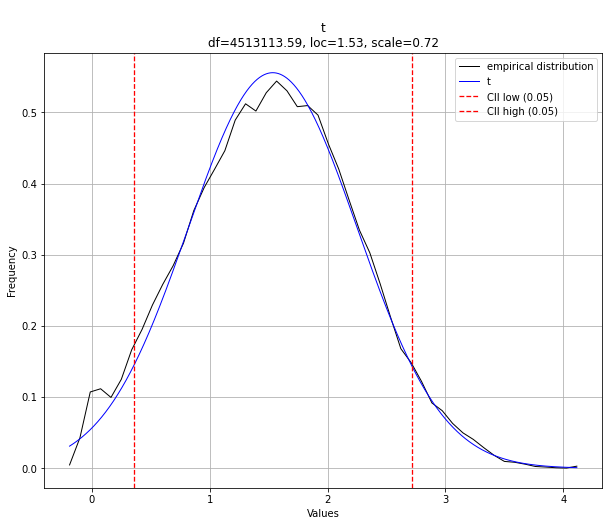}

\caption{Best-fitted probability distributions to the actual data for the two-layer case. (From left to right) First row: $d=5$ and $d=10$; Second row: $d=20$ and $d=50$; Third row: $d=100$ and $d=200$; Fourth row: $d=500$ and $d=1,000$.}
\label{fig:estimated_dist_l2}
\end{figure}

In Tables \ref{tab:estimated_dist_l1} and \ref{tab:estimated_dist_l2}, we present the top 3 distributions for network case, along with the corresponding RSS, while in Figures \ref{fig:estimated_dist} and \ref{fig:estimated_dist_l2}, we depict how close the best-fitted distribution is to the actual data distribution.
\begin{table}[ht]
    \centering
        \caption{Top three fitted probability distributions for the one-layer case. In parentheses, we observe the corresponding RSS.}
    \label{tab:estimated_dist_l1}
    
    \begin{tabular}{cccc}
    \toprule 
        \textbf{Dimension $d$} & \textbf{$1^{st}$} & \textbf{$2^{nd}$} & \textbf{$3^{rd}$}\\
        
        \midrule
    
    5   & Beta ($5.30$) & F ($6.13$)  & Student's t ($6.14$)  \\
    10   & Beta ($0.12$)  & GenExtreme ($0.35$)  & Student's t ($0.37$)  \\
    20   & Student's t  ($0.05$)  & F ($0.8)$& Beta ($0.10$) \\
    50   &  Beta ($0.04$) & Student's t  ($0.20$)  & F ($0.20$)  \\
    100   &  Beta ($0.05$) & Student's t ($0.40$)  & F ($0.41$)  \\
    200   &  Beta ($0.02$) & Student's t ($0.13$)  & F ($0.14$)  \\
    500   &  Beta ($0.01$) & Student's t ($0.06$)  & F ($0.12$)  \\
    1,000   &  Beta ($0.41$) & Student's t ($0.58$)  & F ($0.64$) \\
        \bottomrule
    \end{tabular} 
\end{table}

\begin{table}[ht]
    \centering
        \caption{Top three fitted probability distributions for the two-layer case. In parentheses, we observe the corresponding RSS.}
    \label{tab:estimated_dist_l2}
    
    \begin{tabular}{cccc}
    \toprule 
        \textbf{Dimension $d$} & \textbf{$1^{st}$} & \textbf{$2^{nd}$} & \textbf{$3^{rd}$}\\
        
        \midrule
    
    5   & Beta ($4.61$)   & Student's t ($10.45$) & Exponential ($789.92.23$) \\
    10   & Beta ($0.31$)   & Student's t ($0.77$) & F ($0.78$)   \\
    20   & Beta ($0.02$)   & Student's t ($1.27$) & F ($1.39$)  \\
    50   &  Student's t ($0.04$) &  Beta ($0.05$) & GenExtreme ($1.62$)  \\
    100   & Beta ($0.01$)   & Student's t ($0.02$) & F ($0.02$)  \\
    200   & Beta ($0.004$)   & Student's t ($0.04$) & F ($0.04$) \\
    500   &  Beta ($0.002$)   & Student's t ($0.02$) & F ($0.02$)\\
    1,000   &   Beta ($0.01$)   & Student's t ($0.02$) & F ($0.02$) \\
        \bottomrule
    \end{tabular} 
\end{table}

\par In Tables \ref{tab:estimated_dist_l1} and \ref{tab:estimated_dist_l2}, we observe that the best-fitted distributions are the beta, Student's t, F, and the generalized extreme value. All these distributions satisfy von Mises' condition (\cite{matthys2003estimating}). Furthermore, as shown in Figures \ref{fig:estimated_dist} and \ref{fig:estimated_dist_l2}, the best-fitted distribution is always very close to the empirical distribution obtained from the actual data. This is also verified by the low RSS achieved by the best-fitted distribution in most of the cases. Therefore, we can claim that in all cases the underlying distribution can be described accurately by a distribution that satisfies von Mises' condition.

\par To estimate the extreme value index $\gamma$, we can use the extreme value index closed-form formula for the best-fitted distribution (\cite{de2006extreme}). All of the previously mentioned distributions have a closed-form for calculating the extreme value index. For instance, for beta distribution $\mathcal{B}(\alpha,\beta)$, which is the distribution with the lowest RSS for the majority of the experiments, the closed-form extreme value index is equal to $-\frac{1}{\alpha}$. Since we are optimizing objective functions from ReLU neural networks over a bounded domain, the optimal objective value is always finite. Therefore, intuitively, the beta distribution is a natural choice for the underlying distribution as it has a finite right endpoint. Furthermore, since we are performing uniform sampling and the extreme value index is a function of the parameters of the distribution, we can also create confidence intervals (\cite{ci1,ci2}) around the estimation of the extreme value index.

\end{document}